\documentclass{elsart}
\usepackage{graphics}
\usepackage{graphicx}
\usepackage{epsfig}
\usepackage{amssymb,amsmath}
\usepackage[colorlinks]{hyperref}

\begin{document}
\begin{frontmatter}
\title{New type of monogenic polynomials and associated spheroidal wavelets}
\author{Sabrine Arfaoui}
\address{Department of Informatics, Higher Institute of Applied Sciences and Technology of Mateur, Street of Tabarka, 7030 Mateur, Tunisia.\\
and\\
Research Unit of Algebra, Number Theory and Nonlinear Analysis UR11ES50, Faculty of Sciences, Monastir 5000, Tunisia.}
\ead{arfaoui.sabrine@issatm.rnu.tn}
\author{Anouar Ben Mabrouk}
\address{Higher Institute of Applied Mathematics and Informatics, University of Kairouan, Street of Assad Ibn Alfourat, Kairouan 3100, Tunisia.\\
and\\
Research Unit of Algebra, Number Theory and Nonlinear Analysis UR11ES50, Faculty of Sciences, Monastir 5000, Tunisia.}
\ead{anouar.benmabrouk@fsm.rnu.tn}
\begin{abstract}
In the present work, new classes of wavelet functions are presented in the framework of Clifford analysis. Firstly, some classes of new monogenic polynomials are provided based on 2-parameters weight functions. Such classes englobe the well known Jacobi, Gegenbauer ones. The discovered polynomial sets are next applied to introduce new wavelet functions. Reconstruction formula as well as Fourier-Plancherel rules have been proved.
\end{abstract}
\begin{keyword}
Continuous Wavelet Transform, Clifford analysis, Clifford Fourier transform, Fourier-Plancherel, Monogenic functions.
\PACS: 42B10, 44A15, 30G35.
\end{keyword}
\end{frontmatter}
\section{Introduction}
\qquad Fourier analysis has been for many decades the essential mathematical tool in harmonic analysis and related applications' fields such as physics, engineering, signal/image processing, ... etc. Next, new extending mathematical tool has been introduced to generalize Fourier one and to overcome in some ways the disadvantages of Fourier analysis. It consists of wavelet analysis.

Compared to Fourier theory, wavelets are mathematical functions permitting themselves to cut up data into different components relatively to the frequency spectrum and next focus on these components somehow independently, extract their characteristics and lift to the original data. One main advantage for wavelets is the fact that they are able more than Fourier modes in analyzing discontinuities and/or singularities efficiently and non-stationarity.

Wavelets were developed independently in the fields of mathematics, physics, electrical engineering, and seismic geology. Next, interchanges between these fields have yielded more understanding of their theory and applications.

Nowadays, wavelets are reputable and successful tools in quasi all domains. The particularity in a wavelet basis is the fact that all the elements of a basis are deduced from one source function known as the wavelet mother. Next, such a mother gives raise to all the elements necessary to analyze objects by simple actions of translation, dilatation and rotation. The last parameter is firstly introduced in \cite{Antoine-Murenzi-Vandergheynst} (see also \cite{Antoine-Vandergheynst}) to obtain some directional selectivity of the wavelet transform in higher dimensions and to analyze/characterize spherical data. Indeed, construction of wavelets related to manifolds such as or essentially spheres is based on the geometric structure of the surface where the data lies. This gives raise to the so-called isotropic and anisotropic wavelets.

The present work lies in the whole scope of the study of spherical data. We propose to develop methods based on harmonic structures to define the so-called ultraspheroidal wavelets. One important and actual motivation is issued from 3D-images processing which is noadays a revolutionary task in informatics.

Mathematically, spheroidal functions such as Gegenbauer polynomials which are the starting point in the present extension are solutions with separated variables of the wave equation
$$
\nabla^2w+k^2w=0
$$
in an elliptic cylinder coordinates system, prolate and/or oblate spheroids. In such systems (mainly radial and angular variables), the wave equation above may be transformed to a second order ODE of the form
$$
(1-t^2)w''+2\alpha tw'+(\beta-\gamma^2t^2)w=0.
$$
(See \cite{Abramowitzetal}, \cite{Lietal}, \cite{Morais}, \cite{Osipovetal}, \cite{Stratton}, \cite{Strattonetal}). This last equation leads to special functions such as Bessel, Airy, ... and special polynomials such as Gegenbauer, Legendre, Chebyshev, .... and constitutes a first idea behind the link between these functions and a first motivation of the present work where construction of some new spheroidal mother wavelets are done. Besides, spheroidal functions have been in the basis of modeling physical phenomena where the wave behaviour is pointed out such as radars, antennas, 3D-images, ... Recall also that Gegenbauer polynomials themselves are called ultraspheroidal polynomials. See \cite{Antoine-Murenzi-Vandergheynst}, \cite{Arfaouietal1}, \cite{Arfaouietal2}, \cite{DeSchepper}, \cite{Delanghe}, \cite{Lehar}, \cite{Lietal}, \cite{Michel}, \cite{Moussa}, \cite{Saillardetal}.

The main idea consists in adopting Clifford analysis to introduce or more precisely to extend some existing works on Clifford wavelets for more general cases. Clifford analysis, in its most basic form, is a refinement of harmonic analysis in higher dimensional Euclidean space. By introducing the so-called Dirac operator, researchers introduced the notion of monogenic functions extending holomorphic ones. In this context, different concepts of real and complex analysis have been extended to the Clifford case such as Fourier transform (extended to Clifford Fourier transform, Derivation of functions, ....). For example, Clifford Fourier transform is related or expressed in terms of an exponential operator. For the even dimensional case, it yields a kernel based on Bessel functions. Compared to the classical Fourier transform, the new kenrel satisfies herealso a system of differential equations.

In the present work, one aim is to provide a rigourous development of wavelets adapted to the Clifford calculus. The frame is somehow natural as wavelets are characterized by scale invariance of approximation spaces. Clifford algebra is one mathematical object that owns this characteristic. Recall that multiplication of real numbers scales their magnitudes according to their position in or out from the origin. However, multiplication of the imaginary part of a complex number performs a rotation, it is a multiplication that goes round and round instead of in and out. So, a multiplication of spherical elements by each other results in an element of the sphere. Again, repeated multiplication of the imaginary part results in orthogonal components. Thus, we need a coordinates system that results always in the object, a concept that we will see again and again in the Algebra. In other words, Clifford algebra generalizes to higher dimensions by the same exact principles applied at lower dimensions, by providing an algebraic entity for scalars, vectors, bivectors, trivectors, and there is no limit to the number of dimensions it can be extended to. More details on Clifford algebra, origins, history, developments may be found in \cite{Abreuetal}, \cite{DeBie1}, \cite{DeBie2}, \cite{Delanghe}, \cite{Lehar}, \cite{Pena}.

Let $\Omega$ be an open subset of $\mathbb{R}^m$ or $\mathbb{R}^{m+1}$ and $f:\Omega\rightarrow\mathbb{A}$, where $\mathbb{A}$ is the real Clifford algebra $\mathbb{R}_{m}$ or its complexification $\mathbb{C}_{m}$. $f$ may be written in the form
\begin{equation}\label{fincliffordalbebra}
f=\displaystyle\sum_{A}f_{A}e_{A}
\end{equation}
where the functions $f_A$ are $\mathbb{R}$-valued or $\mathbb{C}$-valued and $(e_A)_A$ is a suitable basis of $\mathbb{A}$.\\

In the literature, there are several techniques available to generate monogenic functions in $\mathbb{R}^{m+1}$ such as the Cauchy-Kowalevski extension (CK-extension) which consists in finding a monogenic extension $g^*$ of an analytic function $g$ defined on a given subset in $\mathbb{R}^{m+1}$ of positive codimension. For analytic functions $g$ on the plane $\{(x_0, \underline{x}) \in\mathbb{R}^{m+1}, \quad x_0 = 0\}$ the problem may be stated as follows: \textit{Find $g^*\in\mathbb{A}$ such that
\begin{equation}\label{cauchy-kowalevski-extension}
\partial_{x_0}g^*=-\partial_{\underline{x}} g^*\quad in \quad \mathbb{R}^{m+1}\quad\hbox{and}\quad
g^*(0, \underline{x}) = g(\underline{x}).
\end{equation}}
A formal solution is
\begin{equation}\label{ck-formal-solution}
g^*(x_0,\underline{x})=\exp(-x_0\partial_{\underline{x}}) g(\underline{x})=\displaystyle\sum_{k=0}^{\infty} \displaystyle\frac{(-x_0)^k}{k!}
\partial_{\underline{x}}^k g(\underline{x}).
\end{equation}

Starting from the real space $\mathbb{R}^m,\;(m>1)$ (or the complex space $\mathbb{C}^m$) endowed with an orthonormal basis $(e_1,\dots, e_m)$, the Clifford algebra $\mathbb{R}_m$ (or its complexifation $\mathbb{C}_m$) starts by introducing a suitable interior product. Let
$$
e_j^2=-1,\quad j=1,\dots,m,
$$
$$
e_je_k+e_ke_j=0,\quad j\neq k,\quad j,k=1,\dots,m.
$$
Two anti-involutions on the Clifford algebra are important. The conjugation is defined as the anti-involution for which
$$
\overline{e_j}=-e_j,\quad j=1,\dots, m
$$
with the additional rule in the complex case,
$$
\overline{i}=-i.
$$
The inversion is defined as the anti-involution for which
$$
e_j^{+}=e_j,\quad j=1,\dots,m.
$$
A basis for the Clifford algebra ($e_A:A\subset\{1,\dots,m\}$)
where $e_{\emptyset}=1$ is the identity element. As these rules are defined, the Euclidian space $\mathbb{R}^m$ is then embedded in the Clifford algebras $\mathbb{R}_m$ and $\mathbb{C}_m$ by identifying the vector $x=(x_1,\dots,x_m)$ with the vector $\underline{x}$ given by
$$
\underline{x}=\displaystyle\sum_{j=1}^{m}e_jx_j.
$$
The product of two vectors is given by
$$
\underline{x}\,\underline{y}=\underline{x}.\underline{y}+\underline{x}\wedge\underline{y}
$$
where
$$
\underline{x}.\underline{y}=-<\underline{x},\underline{y}>=-\displaystyle\sum_{j=1}^{m}x_j\,y_j
$$
and
$$
\underline{x}\wedge\underline{y}=\displaystyle\sum_{j=1}^{m}\displaystyle\sum_{k=j+1}^{m}e_j\,e_k(x_j\,y_k-x_ky_j).
$$
is the wedge product. In particular,
$$
\underline{x}^2=-<\underline{x},\underline{x}>=-|\underline{x}|^2.
$$
An $\mathbb{R}_m$ or $\mathbb{C}_m$-valued function $F(x_1,\dots,x_m)$, respectively $F(x_0, x_1,\dots,x_m)$ is called right monogenic in an open region of $\mathbb{R}^m$, respectively, or $\mathbb{R}^{m+1}$, if in that region
$$
F\partial_{\underline{x}}=0, \quad respectively\quad F(\partial_{x_0}+\partial_{\underline{x}})=0.
$$
Here $\partial_{\underline{x}}$ is the Dirac operator in $\mathbb{R}^m$:
$$
\partial_{\underline{x}}=\displaystyle\sum_{j=1}^{m} e_j \partial_{x_j},
$$
which splits the Laplacian in $\mathbb{R}^m$
$$
\Delta_m=-\partial_{\underline{x}}^2,
$$
whereas $\partial_{x_0}+\partial_{\underline{x}}$ is the Cauchy-Riemann operator in $\mathbb{R}^{m+1}$ for which
$$
\Delta_{m+1}=(\partial_{x_0}+\partial_{\underline{x}})(\partial_{x_0}+\overline{\partial_{\underline{x}}})
$$
Introducing spherical co-ordinates in $\mathbb{R}^m$ by
$$
\underline{x}=r\underline{\omega},\quad r=|\underline{x}|\in[0,+\infty[,\,\underline{\omega}\in S^{m-1},
$$
the Dirac operator takes the form
$$
\partial_{\underline{x}}=\underline{\omega}\left( \partial_r+\displaystyle\frac{1}{r} \Gamma_{\underline{\omega}}\right)
$$
where
$$
\Gamma_{\underline{\omega}}=-\displaystyle\sum_{i<j}e_ie_j(x_i\partial_{x_j}-x_j\partial_{x_i})
$$
is the so-called spherical Dirac operator which depends only on the angular co-ordinates.\\
Throughout this article the Clifford-Fourier transform of $f$ is given by
$$
\mathcal{F}(f(x))(y)=\displaystyle\int_{\mathbb{R}^m} e^{-i<\underline{x},\underline{y}>}\, f(\underline{x}) dV(\underline{x}),
$$
where $dV(\underline{x})$ is the Lebesgues measure on $\mathbb{R}^m$.
\section{A 2-parameters Clifford-Gegenbauer-Jacobi polynomials and associated wavelets}
In this section we propose to introduce a new family of orthogonal polynomials in the Clifford context that generalizes the well-known Jacobi polynomials as well as Clifford-Jacobi polynomials. In the sequel the new polynomials will be denoted by $S_{\ell,m}^{\mu,\alpha}(\underline{x})$. Here, the indexation on $l,m$ is related to the classic indexes relative to the degree and the kind of the polynomial, and $\mu,\alpha$ are related to the new Clifford algebra weight
$$
\omega_{\mu,\alpha}(\underline{x})=|\underline{x}|^{2\mu}(1+|\underline{x}|^2)^\alpha.
$$	
The polynomials $S_{\ell,m}^{\mu,\alpha}(\underline{x})$ are generated as usually by the CK-extension of the monogenicity property of the function $F^*$ which can be written as
\begin{equation}\label{F}
F^{*}(t,\underline{x})= \displaystyle\sum_{\ell=0}^{\infty}\displaystyle\frac{t^\ell}{\ell}S_{\ell,m}^{\mu,\alpha}(\underline{x})\, \omega_{\mu-\ell,\alpha-\ell}(\underline{x}).
\end{equation}
The Dirac operator acting on the CK-extension yields the time derivative of $F^*$ as
$$
\displaystyle\frac{\partial F^*(t,\underline{x})}{\partial t}=\displaystyle\sum_{\ell=0}^{\infty}\displaystyle\frac{t^\ell}{\ell!} S_{\ell+1,m}^{\mu,\alpha}(\underline{x})
\omega_{\mu-\ell-1,\alpha-\ell-1}(\underline{x}),
$$
and the Clifford algebra variable derivative as
$$
\displaystyle\frac{\partial F^*(t,\underline{x})}{\partial \underline{x}}=\displaystyle\sum_{\ell=0}^{\infty}\displaystyle\frac{t^\ell}{\ell!}\left(\partial_{\underline{x}} S_{\ell,m}^{\mu,\alpha}(\underline{x}) \omega_{\mu-\ell,\alpha-\ell}(\underline{x})+ S_{\ell,m}^{\mu,\alpha}(\underline{x}) \partial_{\underline{x}} \omega_{\mu-\ell,\alpha-\ell}(\underline{x}) \right).
$$
Immediate computations yield that
$$
\partial_{\underline{x}} \omega_{\mu-\ell,\alpha-\ell}(\underline{x})
=2\underline{x}\,[(\mu-\ell)\omega_{\mu-\ell-1,\alpha-\ell}(\underline{x}) +(\alpha-\ell) \omega_{\mu-\ell,\alpha-\ell-1}(\underline{x})].
$$
Therefore,
$$
\begin{array}{lll}
\medskip\displaystyle\frac{\partial F^*(t,\underline{x})}{\partial\underline{x}}&=&
\displaystyle\sum_{\ell=0}^{\infty}\,\displaystyle\frac{t^\ell}{\ell!}\left[\partial_{\underline{x}}S_{\ell,m}^{\mu,\alpha}(\underline{x})\, \omega_{\mu-\ell,\alpha-\ell}(\underline{x})\right.\\
\medskip&&\qquad\qquad+S_{\ell,m}^{\mu,\alpha}(\underline{x})
2\underline{x}\,[(\mu-\ell)\omega_{\mu-\ell-1,\alpha-\ell}(\underline{x})\\ \medskip&&\displaystyle\qquad\qquad\qquad\qquad\qquad\left.+(\alpha-\ell)\omega_{\mu-\ell,\alpha-\ell-1}(\underline{x})]\right].
\end{array}
$$
From the monogenicity relation applied to $F$ we derive the recurrence relation
$$
\begin{array}{lll}
\medskip&&S_{\ell+1,m}^{\mu,\alpha}(\underline{x})\omega_{\mu-\ell-1,\alpha-\ell-1}(\underline{x}) +\omega_{\mu-\ell,\alpha-\ell}(\underline{x})\partial_{\underline{x}}S_{\ell,m}^{\mu,\alpha}(\underline{x})\\
\medskip&&+2\underline{x}\,[(\mu-\ell)\omega_{\mu-\ell-1,\alpha-\ell}(\underline{x})
+(\alpha-\ell)\omega_{\mu-\ell,\alpha-\ell-1}(\underline{x})]  S_{\ell,m}^{\mu,\alpha}(\underline{x})=0.
\end{array}
$$
We thus obtain
\begin{equation}\label{14}
\begin{array}{lll}
\medskip S_{\ell+1,m}^{\mu,\alpha}(\underline{x})&=&-2\underline{x}[(\mu-\ell)(1+|\underline{x}|^2)+(\alpha-\ell)|\underline{x}|^2]  S_{\ell,m}^{\mu,\alpha}(\underline{x})\\
\medskip&&-|\underline{x}|^2(1+|\underline{x}|^2)\partial_{\underline{x}}S_{\ell,m}^{\mu,\alpha}(\underline{x}).
\end{array}
\end{equation}
Starting from $S_{0,m}^{\alpha,\mu}(\underline{x})=1$, a straightforward calculation yields for example that
$$
S_{1,m}^{\alpha,\mu}(\underline{x})=-2\underline{x}\,[\mu(1+|\underline{x}|^2)+\alpha|\underline{x}|^2]
=-2[\mu\underline{x}-(\mu+\alpha)\underline{x}^3].
$$
Next, for $\ell=1$ we obtain
$$
\begin{array}{lll}
\medskip S_{2,m}^{\alpha,\mu}(\underline{x})
&=&(4\mu(\mu-1)+2m\mu)\underline{x}^2-[4(\mu-1)(\mu+\alpha)+4\mu(\mu+\alpha-2)\\
\medskip&&+2(\mu+\alpha)(m+2)+2m\mu]\underline{x}^4\\
\medskip&&+[4(\mu+\alpha)(\mu+\alpha-2)+2(\mu+\alpha)(m+2)]\underline{x}^6.
\end{array}
$$
For $\ell=2$ we get
$$
\begin{array}{lll}
\medskip&&S_{3,m}^{\alpha,\mu}(\underline{x})\\
\medskip&=&-2[\mu(2\mu m+4\mu(\mu-1))-2\mu m-8\mu(\mu-1)]\underline{x}^3\\
\medskip&&+[28\mu m+54\mu (\mu-1)+48\mu+12\alpha m+(\mu+\alpha-4)[4\mu m\\
\medskip&&+8\mu(\mu-1)]\underline{x}^5+
[2(\mu-2)[6\mu m+12\mu(\mu-1)+16\mu\alpha\\
\medskip&&+4\alpha m+4\alpha(\alpha-1)]+2(\alpha-2)[4\mu m+8\mu(\mu-1)+8\mu \alpha+2\alpha m]\\
\medskip&&+28\mu m+56\mu(\mu-1)+80\mu+20\alpha m+24\alpha(\alpha-1)
]\underline{x}^7\\
\medskip&&+[2(\mu+\alpha-4)(2\mu m+4\mu(\mu-1)+8\mu\alpha+2\alpha m+4\alpha(\alpha-1))\\
\medskip&&+12\mu m+24\mu(\mu-1)48\mu+12\alpha m+24\alpha(\alpha-1)
]\underline{x}^9.
\end{array}
$$
\begin{rem}
$S_{\ell,m}^{\alpha,\mu}(\underline{x})$ is a polynomial of degree $3\ell$ in $\underline{x}$.
\end{rem}
As for the classical cases, here also we may derive an analogue Rodrigues formulation for the polynomials $S_{\ell,m}^{\mu,\alpha}(\underline{x})$.
\begin{prop}
The $(\mu,\alpha)$ Clifford-Jacobi polynomials $S_{\ell,m}^{\mu,\alpha}(\underline{x})$ may be expressed as
\begin{equation}\label{15}
S_{\ell,m}^{\mu,\alpha}(\underline{x})=(-1)^\ell (|\underline{x}|^2)^{\ell-\mu}(1-\underline{x}^2)^{\ell-\alpha}\, \partial_{\underline{x}}^{\ell} \omega_{\mu,\alpha}(\underline{x}).
\end{equation}
\end{prop}
\hskip-20pt\textbf{Proof.} We proceed by recurrence on $\ell$. For $\ell=0$, we have
$$
\begin{array}{lll}
\medskip|\underline{x}|^{2\mu}  (1-\underline{x}^2)^\alpha&=&
(-1)^0\omega_{\mu-0,\alpha-0}\times 1\\
\medskip&=& (-1)^0\omega_{\mu-0,\alpha-0}\times S_{0,m}^{\alpha,\mu}(\underline{x}).
\end{array}
$$
Thus,
$$
S_{0,m}^{\alpha,\mu}(\underline{x})=(-1)^0\omega_{0-\mu,0-\alpha}(\underline{x})\partial_{\underline{x}}\omega_{\mu,\alpha}(\underline{x}).
$$
For $\ell=1$, we have
$$
\begin{array}{lll}
\medskip\partial_{\underline{x}} (|\underline{x}|^{2\mu}  (1+|\underline{x}|^2)^\alpha)
\medskip&=& 2\mu\underline{x} \omega_{\mu-1,\alpha}(\underline{x})+2\alpha\underline{x}\,\omega_{\mu,\alpha-1}(\underline{x})\\
\medskip&=&(-1)\omega_{\mu-1,\alpha-1}(\underline{x})[-2\mu\underline{x}(1-\underline{x}^2)+2\alpha\underline{x}^3]\\
\medskip&=&(-1)\omega_{\mu-1,\alpha-1}(\underline{x})[-2\mu\underline{x}+2(\mu+\alpha)\underline{x}^3]\\
\medskip&=&(-1)\omega_{\mu-1,\alpha-1}(\underline{x})S_{1,m}^{\mu,\alpha}(\underline{x}).
\end{array}
$$
Thus,
$$
S_{1,m}^{\alpha,\mu}(\underline{x})=(-1)\omega_{1-\mu,1-\alpha}(\underline{x})\partial_{\underline{x}}\omega_{\mu,\alpha}(\underline{x}).
$$
For the convenience, we push the calculus to the order $\ell=2$.
$$
\begin{array}{lll}
\medskip&&\partial_{\underline{x}}^{2}(\,\omega_{\mu,\alpha}(\underline{x}))\\
\medskip&=&2\partial_{\underline{x}}(\mu\underline{x}\omega_{\mu-1,\alpha}(\underline{x})+\alpha\underline{x}\omega_{\mu,\alpha-1}(\underline{x}))\\
\medskip&=&2\mu\,\partial_{\underline{x}}\left( \underline{x}\omega_{\mu-1,\alpha}(\underline{x})\right) +2\alpha\partial_{\underline{x}}(\underline{x} \omega_{\mu,\alpha-1}(\underline{x}))\\
\medskip&=&-2m\mu[\omega_{\mu-1,\alpha}(\underline{x})+4\mu(\mu-1)\underline{x}^2\omega_{\mu-2,\alpha}(\underline{x})+4\alpha \mu \underline{x}^2\omega_{\mu-1,\alpha-1}(\underline{x})\\
\medskip&&-2m\alpha\omega_{\mu,\alpha-1}(\underline{x})+
4\mu\alpha\underline{x}^2\omega_{\mu-1,\alpha-1}(\underline{x})
+4\alpha(\alpha-1)\underline{x}^2\,\omega_{\mu,\alpha-2}(\underline{x})\\
\medskip&=&(-1)^2\omega_{\mu-2,\alpha-2}(\underline{x})[-2m\mu\omega_{1,2}(\underline{x})+4\mu(\mu-1)\underline{x}^2\omega_{0,2}\\
\medskip&&+4\mu\alpha\underline{x}^2\omega_{1,1}(\underline{x})-2m\alpha\omega_{2,1}(\underline{x})+4\mu\alpha\underline{x}^2\omega_{1,1}(\underline{x})\\
\medskip&&+4\alpha(\alpha-1)\underline{x}^2\omega_{2,0}(\underline{x})]\\
\medskip&=& (-1)^2\omega_{\mu-2,\alpha-2}(\underline{x})([2\mu m+4\mu(\mu-1)]\underline{x}^2-[4\mu m+8\mu(\mu-1)\\
\medskip&&+8\mu+2\alpha m]\underline{x}^4+[2\mu m+4\mu(\mu-1)+8\mu+2\alpha m+4\alpha(\alpha-1)]\underline{x}^6)\\
\medskip&=&(-1)^2\omega_{\mu-2,\alpha-2}(\underline{x})S_{2,m}^{\alpha,\mu}(\underline{x}).
\end{array}
$$
Thus,
$$
\begin{array}{lll}
S_{2,m}^{\alpha,\mu}(\underline{x})&=& (-1)^2\omega_{2-\mu,2-\alpha}(\underline{x})\partial_{\underline{x}}^{2}\omega_{\mu,\alpha}(\underline{x}).
\end{array}
$$
Now, assume that
$$
S_{\ell,m}^{\mu,\alpha}(\underline{x})=(-1)^\ell\omega_{\ell-\mu,\ell-\alpha}(\underline{x})\partial_{\underline{x}}^{\ell}
(\omega_{\mu,\alpha}(\underline{x})).
$$
and denote
$$
\Upsilon=-2\underline{x}[(\mu-\ell)(1-\underline{x}^2)+(\alpha-\ell)|\underline{x}|^2] (-1)^\ell\omega_{\ell-\mu,\ell-\alpha}(\underline{x}) \partial_{\underline{x}}^{\ell}\omega_{\mu,\alpha}(\underline{x})
$$
and
$$
\Psi=(-1)^\ell\,\omega_{1,1}(\underline{x})[2(\ell-\mu)\underline{x}\,\omega_{\ell-\mu-1,\ell-\alpha}(\underline{x})
+2m(\ell-\alpha)\underline{x}\omega_{\ell-\mu,\ell-\alpha-1}(\underline{x})].
$$
From equations (\ref{14}) and (\ref{15}) we get
$$
\begin{array}{lll}
\medskip&&S_{\ell+1,m}^{\mu,\alpha}(\underline{x})\\
\medskip&=&-2\underline{x}[(\mu-\ell)(1-\underline{x}^2)+(\alpha-\ell)|\underline{x}|^2] (-1)^\ell\omega_{\ell-\mu,\ell-\alpha}(\underline{x}) \partial_{\underline{x}}^{\ell}\omega_{\mu,\alpha}(\underline{x})\\
\medskip&&-\omega_{1,1}(\underline{x})\partial_{\underline{x}}((-1)^\ell\omega_{\ell-\mu,\ell-\alpha}(\underline{x})
\partial_{\underline{x}}^{\ell}\omega_{\mu,\alpha}(\underline{x}))\\
\medskip&=&\Upsilon-(-1)^\ell\omega_{1,1}(\underline{x})\partial_{\underline{x}}[\omega_{\ell-\mu,\ell-\alpha}(\underline{x}) \partial_{\underline{x}}^{\ell}\omega_{\mu,\alpha}(\underline{x})]\\
\medskip&=&\Upsilon-\Psi\partial_{\underline{x}}^{\ell}\omega_{\mu,\alpha}(\underline{x})-(-1)^\ell\omega_{\ell-\mu+1,\ell-\alpha+1}(\underline{x}) \partial_{\underline{x}}^{\ell+1}\omega_{\mu,\alpha}(\underline{x}).
\end{array}
$$
On the other hand, we have
$$
\begin{array}{lll}
\medskip\Psi&=&(-1)^\ell\omega_{1,1}(\underline{x})[2(\ell-\mu)\underline{x}\omega_{\ell-\mu-1,\ell-\alpha}(\underline{x})+2(\ell-\alpha)\underline{x}\omega_{\ell-\mu,\ell-\alpha-1}(\underline{x})]\\
\medskip&=&(-2\underline{x})(-1)^\ell[(\mu-\ell)(1-\underline{x}^2)-(\alpha-\ell)\underline{x}^2]\omega_{\ell-\mu,\ell-\alpha}(\underline{x})\\
\medskip&=&\displaystyle\frac{\Upsilon}{\partial_{\underline{x}}^{\ell}\omega_{\mu,\alpha}(\underline{x})}.
\end{array}
$$
Hence, we get
$$
S_{\ell+1,m}^{\mu,\alpha}(\underline{x})=(-1)^{\ell+1}\omega_{\ell-\mu+1,\ell-\alpha+1}(\underline{x}) \partial_{\underline{x}}^{\ell+1}\omega_{\mu,\alpha}(\underline{x}).
$$
\begin{prop}\label{Orthogonalityrelation1} Let the integral
$$
I_{\ell,t,p}^{\mu,\alpha}=\displaystyle\int_{\mathbb{R}^m}\underline{x}^{\ell}S_{t,m}^{\mu+p,\alpha+p}(\underline{x})\,\omega_{\mu,\alpha}(\underline{x})\,dV(\underline{x}).
$$
Then, the following orthogonality relation holds.
\begin{equation}\label{16}
I_{\ell,t,t}^{\mu,\alpha}=0
\end{equation}	
whenever $4t<1-m-2(\mu+\alpha)$.
\end{prop}
\hskip-20pt\textbf{Proof.} Denote
$$
I_{\ell,t}=\displaystyle\int_{\mathbb{R}^m}\underline{x}^{\ell}\partial_{\underline{x}}^t(\omega_{\mu+t,\alpha+t}(\underline{x}))dV(\underline{x}).
$$
Using Stokes's theorem, we obtain
$$
\begin{array}{lll}
\medskip&&\displaystyle\int_{\mathbb{R}^m}\underline{x}^{\ell} S_{t,m}^{\mu+t,\alpha+t}(\underline{x})\omega_{\mu,\alpha}(\underline{x})dV(\underline{x})\\ \medskip&=&\displaystyle\int_{\mathbb{R}}\underline{x}^{\ell}(-1)^t\omega_{t-\mu-t,t-\alpha-t}(\underline{x})
\partial_{\underline{x}}^t(\omega_{\mu+t,\alpha+t}(\underline{x}))\omega_{\mu,\alpha}(\underline{x})dV(\underline{x})\\
\medskip&=&(-1)^t\displaystyle\int_{\mathbb{R}^m}\underline{x}^{\ell}\partial_{\underline{x}}^t(\omega_{\mu+t,\alpha+t}(\underline{x}))dV(\underline{x})\\
\medskip&=&(-1)^t\displaystyle\int_{\mathbb{R}^m}\underline{x}^{\ell}\partial_{\underline{x}}\partial_{\underline{x}}^{t-1}(\omega_{\mu+t,\alpha+t}(\underline{x}))dV(\underline{x})\\
\medskip&=&(-1)^t\left[\displaystyle\int_{\partial\mathbb{R}^m}\underline{x}^{\ell}\partial_{\underline{x}}^{t-1}(\omega_{\mu+t,\alpha+t}(\underline{x}))\partial\Gamma(\underline{x})\right.\\
\medskip&&\qquad\qquad\qquad\left.-\displaystyle\int_{\mathbb{R}^m}\partial_{\underline{x}}(\underline{x}^{\ell}) \partial_{\underline{x}}^{t-1}(\omega_{\mu+t,\alpha+t}(\underline{x}))dV(\underline{x})\right].
\end{array}
$$
Denote
$$
I=\displaystyle\int_{\partial\mathbb{R}^m}\underline{x}^{\ell}\partial_{\underline{x}}^{t-1}(\omega_{\mu+t,\alpha+t}
(\underline{x}))\partial\Gamma(\underline{x})
$$
and
$$
II=\displaystyle\int_{\mathbb{R}^m}\partial_{\underline{x}}(\underline{x}^{\ell})\partial_{\underline{x}}^{t-1}(\omega_{\mu+t,\alpha+t}(\underline{x}))dV(\underline{x}).
$$
The integral $I$ vanishes due to the asumption  $4t<1-m-2(\mu+\alpha)$. We now apply the following technical result.
\begin{lem}\label{lemmaevenodd}
For all $n\in\mathbb{N}$, we have
\begin{equation}
\partial_{\underline{x}}(\underline{x}^n)=\gamma_{n,m} \underline{x}^{n-1},
\end{equation}
where
$$
\gamma_{n,m}=
\begin{cases}
-n\quad\hbox{if}\quad n\,\mbox{is \,even},\\
-(m+n-1)\quad\mbox{if}\quad n\,\mbox{is\, odd}.
\end{cases}
$$
\end{lem}
\hskip-20pt\,Due to Lemma \ref{lemmaevenodd}, we get
$$
\begin{array}{lll}
II&=&\gamma_{l,m}\displaystyle\int_{\mathbb{R}^m}\underline{x}^{\ell-1}\partial_{\underline{x}}^{t-1}(\omega_{\mu+t,\alpha+t}(\underline{x}))dV(\underline{x})\\
&=&\gamma_{l,m}I_{\ell-1,t-1}.
\end{array}
$$
Hence we obtain
$$
\begin{array}{lll}
\medskip&&\displaystyle\int_{\mathbb{R}^m}\underline{x}^{\ell} S_{t,m}^{\mu+t,\alpha+t}(\underline{x})\omega_{\mu,\alpha}(\underline{x})dV(\underline{x})\\
\medskip&=&(-1)^{t+1}\gamma_{l,m} I_{\ell-1,t-1}\\
\medskip&=&(-1)^{t+1}\gamma_{l,m}[(-1)^{t}\gamma_{l-1,m}I_{\ell-2,t-2}]\\
\medskip&=&(-1)^{2t+1}\gamma_{l,m}\gamma_{l-1,m}I_{\ell-2,t-2}\\
&\vdots&\\
&=& C(m,\ell,t) I_0\\
&=&0,
\end{array}
$$
where
\begin{equation}\label{C(m,l,t)}
C(m,\ell,t)=(-1)^{ml+1}\displaystyle\prod_{k=0}^{m}\gamma_{k,m}.
\end{equation}
We now introduce the generalized $(\mu,\alpha)$-Clifford-Gegenbauer-Jacobi wavelets associated to the polynomials introduced previously. Note that proposition \ref{Orthogonalityrelation1} implies that for $0<t<\displaystyle\frac{1-m-2(\mu+\alpha)}{4}$,
$$
\displaystyle\int_{\mathbb{R}^m} S_{t,m}^{\mu+t,\alpha+t}(\underline{x})\omega_{\alpha,\beta}(\underline{x})dV(\underline{x})=0.
$$
\begin{defn} The generalized $(\mu,\alpha)$ Clifford-Gegenbauer-Jacobi analyzing wavelet is defined by
$$
\psi_{\ell,m}^{\mu,\alpha}(\underline{x})=S_{\ell,m}^{\mu+\ell,\alpha+\ell}(\underline{x}) \,\omega_{\mu,\alpha}(\underline{x})
=(-1)^\ell \partial_{\underline{x}}^{\ell} \omega_{\mu+\ell,\alpha+\ell}(\underline{x}).
$$
\end{defn}
The wavelet $\psi_{\ell,m}^{\mu,\alpha}$ has vanishing moments as is shown in the next proposition.
\begin{prop} The following assertions are true.
\begin{enumerate}
\item Whenever $0<k<-m-\ell-2(\mu+\alpha)$ and $k<\ell$ we have
\begin{equation}\label{17}
\displaystyle\int_{\mathbb{R}^m} \underline{x}^k \psi_{\ell,m}^{\mu,\alpha}(\underline{x}) dV(\underline{x})=0.
\end{equation}
\item The Clifford-Fourier transform of $\psi_{\ell,m}^{\mu,\alpha}$ takes the form
\begin{equation}
\widehat{\psi_{\ell,m}^{\mu,\alpha}(\underline{u})}=(-i)^\ell\, \xi^\ell (2\pi)^{\frac{m}{2}}\rho^{1-\frac{m}{2}+\ell}\,\displaystyle\int_{0}^\infty\widetilde{\omega}^{l,m}_{\alpha,\mu}(r)\,J_{\frac{m}{2}-1}(r\rho)dr
\end{equation}
where
$$
\widetilde{\omega}^{l,m}_{\alpha,\mu}(r)=r^{2(\mu+\ell)+\frac{m}{2}} (1+r^2)^{\alpha+\ell}.
$$
\end{enumerate}
\end{prop}
\hskip-20pt\textbf{Proof.} The first assertion is a natural consequence of proposition \ref{Orthogonalityrelation1}. We prove the second. We have
$$
\widehat{\psi}_{\ell,m}^{\mu,\alpha}(\underline{u})=(-1)^\ell (i\underline{u})^\ell \widehat{\omega_{\mu+\ell,\alpha+\ell}}(\underline{u}).
$$
This Fourier transform  can be simplified by using the spherical co-ordinates. By definition, we have
\begin{equation}\label{spherical-co-ordinates}
\widehat{\omega_{\mu+\ell,\alpha+\ell}}(\underline{u})=\displaystyle\int_{\mathbb{R}^m} |\underline{x}|^{2(\mu+\ell)} (1+|\underline{x}|^2)^{\alpha+\ell}\, e^{-i<\underline{x},\underline{u}>} dV(\underline{x}).
\end{equation}
Introducing spherical co-ordinates
$$
\underline{x}=r\underline{\omega},\quad \underline{u}=\rho\underline{\xi},\quad r=|\underline{x}|,\quad \rho=|\underline{u}|, \quad \underline{\omega}\in S^{m-1}, \,\underline{\xi}\in S^{m-1}
$$
(where $S^{m-1}$ is the unit sphere of $\mathbb{R}^m$) expression (\ref{spherical-co-ordinates}) becomes
$$
\begin{array}{lll}
\widehat{\omega_{\mu+\ell,\alpha+\ell}(\underline{u})}&=&\displaystyle\int_{0}^\infty r^{2(\mu+\ell)+m-1} (1+r^2)^{\alpha+\ell} \,dr\displaystyle\int_{S^{m-1}} e^{-i<r\underline{\omega},\rho\underline{\xi}>} d\sigma(\underline{\omega})
\end{array}
$$
where $d\sigma(\underline{\omega})$ stands for the Lebesgue measure on $S^{m-1}$.\\
We now use the following technical result which is known in the theory of Fourier analysis of radial functions and the theory of bessel functions.
\begin{lem}\label{BesselFourierTransform}\cite{Stein-Weiss}
$$
\displaystyle\int_{S^{m-1}}  e^{-i<r\underline{\omega},\rho\underline{\xi}>} d\sigma(\underline{\omega})=\displaystyle\frac{(2\pi)^{\frac{m}{2}} J_{\frac{m}{2}-1}(r\rho)}{(r\rho)^{\frac{m}{2}-1}}
$$
where $J_{\frac{m}{2}-1}$ is the bessel function of the first kind of order $\frac{m}{2}-1$ and $d\sigma$ is the Lebesgue measure on the sphere $S^{m-1}$.
\end{lem}
\hskip-20pt Now, according to lemma \ref{BesselFourierTransform}, we obtain
$$
\begin{array}{lll}
\medskip\widehat{\omega_{\mu+\ell,\alpha+\ell}(\underline{u})}
&=&(2\pi)^{\frac{m}{2}}\rho^{1-\frac{m}{2}}\,\displaystyle\int_{0}^\infty\,r^{2(\mu+\ell)+\frac{m}{2}}(1+r^2)^{\alpha+\ell}\,J_{\frac{m}{2}-1}(r\rho)dr\\
\medskip&=&(2\pi)^{\frac{m}{2}}\rho^{1-\frac{m}{2}}\,\displaystyle\int_{0}^\infty\widetilde{\omega}^{l,m}_{\alpha,\mu}(r)\,J_{\frac{m}{2}-1}(r\rho)dr.
\end{array}
$$
Consequently, we obtain the following expression for the Fourier transform of the $(\mu,\alpha)$-clifford-jacobi wavelets
$$
\widehat{\psi_{\ell,m}^{\mu,\alpha}}(\underline{u})=(-i)^\ell\,\xi^\ell(2\pi)^{\frac{m}{2}}\rho^{1-\frac{m}{2}+\ell}\,\displaystyle\int_{0}^\infty\widetilde{\omega}^{l,m}_{\alpha,\mu}(r)\,J_{\frac{m}{2}-1}(r\rho)dr.
$$
\textbf{Proof of Lemma \ref{BesselFourierTransform}.}
Denote firstly $I$ the right hand side integral. Recall next that for $\lambda\in\mathbb{R}$, the Bessel function $J_{\lambda}$ may be written in the integral form
$$
J_{\lambda}(x)=\dfrac{1}{2\pi}\displaystyle\int_{-\pi}^{\pi}e^{-i(\lambda t-xsint)}dt.
$$
Next, as the mesure $d\sigma$ is invariant under rotations, we may assume without loss of the generality that $\xi=(1,0,\dots,0)$. As a result, we obtain
$$
I(\xi)=\displaystyle\int_{S^{m-1}}e^{-ir\rho\cos\theta}d\sigma(\omega),
$$
where $\theta$ is the angle $(\omega,e_1)$, with $e_1=(1,0,\dots,0)$. In spherical coordinates, this means that
$$
I(\xi)=\omega_{m-1}\displaystyle\int_{0}^{\pi}e^{-ir\rho\cos\theta}\sin^{m-2}\theta\,d\theta.
$$
Denote next $t=\cos\theta$. We get
$$
I(\xi)=\omega_{m-1}\displaystyle\int_{-1}^{1}e^{-ir\rho t}(1-t^2)^{(m-3)/2}dt.
$$
Observing next that the area $\omega_{m-1}$ of the unit sphere $S^{m-1}$ is
$$
\omega_{m-1}=\dfrac{\pi^{\frac{m-3}{2}}}{\Gamma(\frac{m-1}{2})}
$$
and in the other hand,
$$
J_\lambda(x)=\dfrac{|\frac{x}{2}|^{\lambda}}{\pi^{1/2}\Gamma(\lambda+1/2)}\displaystyle\int_{-1}^{1}e^{-ixt}(1-t^2)^{\lambda-1/2}dt,
$$
we obtain the desired result.

A first question in wavelet theory is the admissibility of the wavelet mother. This is checked in the following lemma.
\begin{lem}
The wavelet mother $\psi_{\ell,m}^{\mu,\alpha}$ satisfies the dmissibility assumption
\begin{equation}\label{41}
{\mathcal{A}_{\ell,m}^{\mu,\alpha}}=\displaystyle\frac{1}{\omega_{m-1}}\displaystyle\int_{\mathbb{R}^m}\left|\widehat{\psi_{\ell,m}^{\mu,\alpha}}(\underline{x})\right|^2 \displaystyle\frac{dV(\underline{x})}{|\underline{x}|^m}<+\infty.
\end{equation}
\end{lem}
\hskip-20pt Indeed,

Now, we introduce the Generalized $(\mu,\alpha)$-Clifford-Gegenbauer-Jacobi Continuous Wavelet Transform. For $a>0$ and $\underline{b}\in\mathbb{R}_m$, the $(a,\underline{b})$-copy of the wavelet mother $\psi_{\ell,m}^{\mu,\alpha}$ is defined by
\begin{equation} \label{11}
^a_{\underline{b}}{\psi}_{\ell,m}^{\mu,\alpha}(\underline{x})=a^{-\frac{m}{2}}{\psi}_{\ell,m}^{\mu,\alpha}(\displaystyle\frac{\underline{x}-\underline{b}}{a}).
\end{equation}
\begin{defn}\label{Generalized-mu-alpha-Clifford-Jacobi-Continuous-Wavelet-Transform}
The generalized $(\mu,\alpha)$-Clifford-Gegenbauer-Jacobi CWT  applies to functions $f\in L_2$ by means of the wavelet coefficient
$$
C_{a,\underline{b}}(f)= <^a_{\underline{b}}\!\!{\psi}_{\ell,m}^{\mu,\alpha},f>
=\displaystyle\int_{\mathbb{R}^m}f(\underline{x})\,
^a_{\underline{b}}{\psi}_{\ell,m}^{\mu,\alpha}(\underline{x})dV(\underline{x}).
$$
\end{defn}
\hskip-20pt Introducing the inner product
$$
<C_{a,\underline{b}}(f),C_{a,\underline{b}}(g)>
=\displaystyle\frac{1}{\mathcal{A}_{\ell,m}^{\mu,\alpha}}\displaystyle\int_{\mathbb{R}^m}\displaystyle\int_{0}^{+\infty}
\overline{C_{a,\underline{b}}(f)}C_{a,\underline{b}}(g)\displaystyle\frac{da}{a^{m+1}}dV(\underline{b})
$$
we obtain the following result.
\begin{thm}\label{ReconstructionFormula1}
Any function $f$ in $L_2(\mathbb{R}_m)$ may be reconstructed in the $L_2$-sense as
\begin{equation}\label{13}
f(x)=\displaystyle\frac{1}{{\mathcal{A}_{\ell,m}^{\mu,\alpha}}}\displaystyle\int_{a>0}\displaystyle\int_{b\in\mathbb{R}^m}C_{a,\underline{b}}(f)\psi\left(\displaystyle\frac{\underline{x}-\underline{b}}{a}\right)\displaystyle\frac{da\,dV(\underline{b})}{a^{m+1}}.
\end{equation}
\end{thm}
The proof reposes on the following result.
\begin{lem}\label{FourierPlancherel1}
It holds that
$$
\displaystyle\int_{a>0}\displaystyle\int_{b\in\mathbb{R}^m}\overline{C_{a,\underline{b}}(f)}C_{a,\underline{b}}(g)\displaystyle\frac{da\,dV(\underline{b})}{a^{m+1}}
=\mathcal{A}_{\ell,m}^{\mu,\alpha}\displaystyle\int_{\mathbb{R}^m}f(\underline{x})\overline{g(\underline{x})}dV(\underline{x}).
$$
\end{lem}
\hskip-20pt\textbf{Proof.} Using the Clifford Fourier transform we observe that
$$
C_{a,\underline{b}}(f)(\underline{b})=\widetilde{a^{\frac{m}{2}}\widehat{\widehat{f}(\underline{.})\widehat{\psi}(a\underline{.})}}(\underline{b}),
$$
where, $\widetilde{h}(\underline{u})=h(-\underline{u})$, $\forall\,h$. Thus,
$$
\overline{C_{a,\underline{b}}(f)}C_{a,\underline{b}}(g)=\overline{\widehat{\left(\widehat{f}(\underline{.})a^{\frac{m}{2}} \widehat{\psi}(a\underline{.})\right)}}(-\underline{b}) \widehat{\left(\widehat{g}(\underline{.})a^{\frac{m}{2}} \widehat{\psi}(a\underline{.})\right)}(-\underline{b}).
$$
Consequently,
$$
\begin{array}{lll}
\medskip<C_{a,\underline{b}}(f),C_{a,\underline{b}}(g)>
&=&\displaystyle\int_{a>0}\displaystyle\int_{\mathbb{R}^m}\overline{\widehat{\widehat{f}(\underline{.})a^{\frac{m}{2}}
\widehat{\psi}(a\underline{.})}}\,\widehat{\widehat{g}(\underline{.})a^{\frac{m}{2}}\widehat{\psi}(a\underline{.})}\displaystyle\frac{da\,dV(\underline{b})}{a^{m+1}}\\
\medskip&=&\displaystyle\int_{a>0}\displaystyle\int_{\mathbb{R}^m} \overline{\widehat{f}(\underline{b})}\widehat{g}(\underline{b})
\displaystyle\frac{a^m|\widehat{\psi}(a\underline{b})|^2}{a^{m+1}}\,da\,dV(\underline{b})\\
\medskip&=&{\mathcal{A}_{\ell,m}^{\mu,\alpha}}\displaystyle\int_{\mathbb{R}^m}\overline{\widehat{f(\underline{b})}}\widehat{g}(b)dV(\underline{b})\\
\medskip&=&{\mathcal{A}_{\ell,m}^{\mu,\alpha}}<\widehat{f},\widehat{g}>\\
&=&<f,g>.
\end{array}
$$
\textbf{Proof of Theorem \ref{ReconstructionFormula1}.} It follows immediately from lemma \ref{FourierPlancherel1} and Riesz rule.
\section{A 2-parameters Clifford-Gauss-Gegenbauer-Jacobi Polynomials and associated wavelets}
In this section, we develop a second class of orthogonal polynomials and associated wavelets already generalising the well known classes of Gegenbauer and Jacobi and also based on a 2-parameters weight function on the Clifford algebra. Polynomials elements of the new class will be denoted by $K_{\ell,m}^{\alpha,-\beta}(\underline{x})$. These are generated by the CK-extension $F^*$ of the 2-parameters weight function
$$
\omega_{\alpha,-\beta}(\underline{x})=(1+|\underline{x}|^2)^\alpha e^{-\beta|\underline{x}|^2},\quad \beta>0.
$$
The CK-extension then takes the form
$$
\begin{array}{lll}
F^{*}(t,\underline{x})
&=& \displaystyle \sum_{\ell=0}^{\infty}\displaystyle\frac{t^\ell}{\ell!}K_{\ell,m}^{\alpha,-\beta}(\underline{x})\, \omega_{\alpha-\ell,-\beta}(\underline{x}).
\end{array}
$$
The Dirac operator acting on the CK-extension can be written as
$$
\displaystyle\frac{\partial F^*(t,\underline{x})}{\partial t}=\displaystyle\sum_{\ell=0}^{\infty}\displaystyle\frac{t^\ell}{\ell!} K_{\ell+1,m}^{\alpha,-\beta}(\underline{x})\, \omega_{\alpha-\ell-1, -\beta}(\underline{x})
$$
and
$$
\begin{array}{lll}
\medskip\displaystyle\frac{\partial F^*(t,\underline{x})}{\partial \underline{x}}&=&\displaystyle\sum_{\ell=0}^{\infty}\displaystyle\frac{t^\ell}{\ell!}\left[\omega_{\alpha-\ell,-\beta}(\underline{x})\partial_{\underline{x}}K_{\ell,m}^{\alpha,-\beta}(\underline{x})\right.\\
\medskip&&\qquad\qquad+[2(\alpha-\ell)\underline{x}\, \omega_{\alpha-\ell-1, -\beta}(\underline{x})\\
\medskip&&\qquad\qquad\left.-2\beta\underline{x}\, \omega_{\alpha-\ell, -\beta}(\underline{x})]\,K_{\ell,m}^{\alpha,-\beta}(\underline{x})\right].
\end{array}
$$
From the monogenicity relation, we get
$$
\begin{array}{lll}
\medskip&&K_{\ell+1,m}^{\alpha,-\beta}(\underline{x})\, \, \omega_{\alpha-\ell-1, -\beta}(\underline{x})+
\, \omega_{\alpha-\ell, -\beta}(\underline{x})\partial_{\underline{x}}K_{\ell,m}^{\alpha,-\beta}(\underline{x})\\
\medskip&&+[2(\alpha-\ell)\underline{x}\, \omega_{\alpha-\ell-1, -\beta}(\underline{x})-
2\beta\underline{x}\, \omega_{\alpha-\ell, -\beta}(\underline{x})]\,K_{\ell,m}^{\alpha,-\beta}(\underline{x})\\
\medskip&&=0.
\end{array}
$$
Finally, the following recurrence relation is obtained.
\begin{equation}\label{19}
\begin{array}{lll}
\medskip K_{\ell+1,m}^{\alpha,-\beta}(\underline{x})
&=&-[2(\alpha-\ell)\underline{x}-2\beta\underline{x}\, \omega_{1,0}(\underline{x})]\,K_{\ell,m}^{\alpha,-\beta}(\underline{x})\\
\medskip&&-\omega_{1,0}(\underline{x})\,\partial_{\underline{x}}K_{\ell,m}^{\alpha,-\beta}(\underline{x}).
\end{array}
\end{equation}
As $K_{0,m}^{1,-\beta}(\underline{x})=1$, a straightforward calculation yields that
$$
K_{1,m}^{\alpha,-\beta}(\underline{x})=-2(\alpha-\beta)\underline{x}-2\beta\underline{x}^3.
$$
Next, for $\ell=1$, we obtain
$$
\begin{array}{lll}
\medskip K_{2,m}^{\alpha,-\beta}(\underline{x})
&=&\partial_{\underline{x}}[2\alpha\underline{x}\,\omega_{\alpha-1,-\beta}(\underline{x})-2\beta\underline{x}\,\omega_{\alpha,-\beta}(\underline{x})]\\
\medskip&=&[4(\alpha-1-\beta)(\alpha-\beta)-2\beta(m+2)+2m(\alpha-\beta)]\underline{x}^2\\
\medskip&+&[4(\alpha-1-\beta)\beta+4(\alpha-\beta)\beta+2\beta(m+2)]\underline{x}^4\\
\medskip&+&4\beta^2\underline{x}^6-2m(\alpha-\beta)\\
\medskip&=&[2\alpha m+4\alpha(\alpha-1)-8\alpha\beta-4m\beta+4\beta^2]\underline{x}^2\\
\medskip&&+(8\alpha\beta+2m\beta-8\beta^2)\underline{x}^4+ 4\beta^2\underline{x}^6-2m(\alpha-\beta).
\end{array}
$$
For $\ell=2$, we get
$$
\begin{array}{lll}
\medskip&&K_{3,m}^{\alpha,-\beta}(\underline{x})\\
\medskip&=&[4m(\alpha^2+\beta^2)+4(2\alpha^2+\beta^2)+4m(\beta-\alpha)-8\alpha-8\alpha\beta(1+m)]\underline{x}\\
\medskip&&-8\alpha\beta(1+m)]\underline{x}
+[-8(\alpha^3-\beta^3)+8\alpha\beta(1+m)+4\alpha m(1-\alpha)\\
\medskip&&-24\beta^2(1+\alpha)-8\alpha(1-2\alpha)+24\alpha^2\beta-12m\beta^2]\underline{x}^3\\
\medskip&&+[-24\alpha^2\beta-8m\alpha\beta+48\alpha\beta^2+12m\beta^2-24\beta^3+8\alpha\beta+24\beta^2)\underline{x}^5\\
\medskip&&+(-24\alpha\beta^2-8\beta^2-8\beta^3-4m\beta^2)\underline{x}^7+8\beta^3\underline{x}^9.
\end{array}
$$
Remark that $K_{\ell,m}^{\alpha,-\beta}(\underline{x})$ is a polynomial of degree $3\ell$ in $\underline{x}$.
\begin{prop}
The Generalized Clifford-Gauss-Gagenbauer-Jacobi polynomials $K_{\ell,m}^{\alpha,-\beta}(\underline{x})$ may be expressed by
\begin{equation}\label{20}
K_{\ell,m}^{\alpha,-\beta}(\underline{x})=(-1)^\ell \omega_{\ell-\alpha,\beta}(\underline{x})\,
\partial_{\underline{x}}^{\ell}(\omega_{\alpha,-\beta}(\underline{x})).
\end{equation} 	
\end{prop}
\hskip-20pt\textbf{Proof.} For $\ell=0$, we have
$$
K_{0,m}^{\alpha,-\beta}(\underline{x})=1
$$
and on the right hand side we have
$$
\omega_{-\alpha,\beta}(\underline{x})\,
\omega_{\alpha,-\beta}(\underline{x})=1.
$$
For $\ell=1$, we have
$$
\begin{array}{lll}
\medskip\partial_{\underline{x}}^{1}(\omega_{\alpha,-\beta}(\underline{x}))&=&2\alpha\underline{x}\,\omega_{\alpha-1,-\beta}(\underline{x})-2\beta\underline{x}\,\omega_{\alpha,-\beta}(\underline{x})\\
\medskip&=&(-1)\,\omega_{\alpha-1,-\beta}(\underline{x})(-2(\alpha-\beta)\underline{x}-2\beta\underline{x}^3)\\
\medskip&=&(-1)\,\omega_{\alpha-1,-\beta}(\underline{x}){K_{1,m}^{\alpha,-\beta}(\underline{x})}.
\end{array}
$$
Therefore,
$$
K_{1,m}^{\alpha,-\beta}(\underline{x})=(-1)\omega_{1-\alpha,\beta}(\underline{x})\,\partial_{\underline{x}}^{1}(\omega_{\alpha,-\beta}(\underline{x})).
$$
Now, as previously, we explain the case $\ell=2$. We have
$$
\begin{array}{lll}
\medskip&&\partial_{\underline{x}}^{2}(\omega_{\alpha,-\beta}(\underline{x}))\\
\medskip&=&\partial_{\underline{x}}[2\alpha\underline{x}\omega_{\alpha-1,-\beta}(\underline{x})-2\beta\underline{x}\omega_{\alpha,-\beta}(\underline{x})]\\
\medskip&=&-2m\alpha\,\omega_{\alpha-1,-\beta}(\underline{x})+4\alpha(\alpha-1)\underline{x}^2\omega_{\alpha-2,-\beta}(\underline{x})\\
\medskip&-&4m^2\alpha\beta\underline{x}^2\omega_{\alpha-1,-\beta}(\underline{x})+ 2m^2\beta\omega_{\alpha,-\beta}(\underline{x})\\
\medskip&-&4\alpha\beta\underline{x}^2\,\omega_{\alpha-1,-\beta}(\underline{x})+4\beta^2\,\omega_{\alpha,-\beta}(\underline{x})\\
\medskip&=&(-1)^2\,\omega_{\alpha-2,-\beta}(\underline{x})
[2\alpha m+4\alpha(\alpha-1)-8\alpha\beta-4m\beta+4\beta^2] \underline{x}^2\\
\medskip&&+(8\alpha\beta+2m\beta-8\beta^2)\underline{x}^4+4\beta^2\underline{x}^6-2\alpha m+2m\beta]\\
\medskip&=&(-1)^2\,\omega_{\alpha-2,-\beta}(\underline{x})\,K_{2,m}^{\alpha,-\beta}(\underline{x}).
\end{array}
$$
Consequently,
$$
K_{2,m}^{\alpha,-\beta}(\underline{x})=(-1)^2\,\omega_{2-\alpha,\beta}(\underline{x}) \,\partial_{\underline{x}}^{2}(\omega_{\alpha,-\beta}(\underline{x})).
$$
Denote
$$
\Lambda(\underline{x})=-[2(\alpha-\ell)\underline{x}-2\beta\underline{x}\omega_{1,0}(\underline{x})](-1)^\ell \omega_{\ell-\alpha,\beta}(\underline{x})\,
\partial_{\underline{x}}^{\ell}{\omega_{\alpha,-\beta}(\underline{x})},
$$
and
$$
\Theta(\underline{x})= (-1)^\ell\omega_{1,0}(\underline{x})[2(\ell-\alpha)\underline{x}\omega_{\ell-\alpha-1,\beta}(\underline{x})+2\beta\underline{x}\omega_{\ell-\alpha,\beta}(\underline{x})]\,\partial_{\underline{x}}^{\ell}(\omega_{\alpha,-\beta}(\underline{x})).
$$
From equations (\ref{19}) and (\ref{20}) we observe that
$$
\begin{array}{lll}
\medskip K_{\ell+1,m}^{\alpha,-\beta}(\underline{x})
&=&-[2(\alpha-\ell)\underline{x}-2\beta\underline{x}\omega_{1,0}(\underline{x})](-1)^\ell \omega_{\ell-\alpha,\beta}(\underline{x})\,
\partial_{\underline{x}}^{\ell}{\omega_{\alpha,-\beta}(\underline{x})}\\
\medskip&\,& \omega_{1,0}(\underline{x})\partial_{\underline{x}}((-1)^\ell \omega_{\ell-\alpha,\beta}(\underline{x})\,
\partial_{\underline{x}}^{\ell}\omega_{\alpha,-\beta}(\underline{x}))\\
\medskip&=&\Lambda(\underline{x})-\Theta(\underline{x}) +(-1)^{\ell+1}\omega_{\ell-\alpha+1,\beta}(\underline{x})\,\partial_{\underline{x}}^{\ell+1}(\omega_{\alpha,-\beta}(\underline{x})).
\end{array}
$$
In the other hand, we have
$$
\begin{array}{lll}
\medskip\Theta(\underline{x})
&=&(-1)^\ell[2(\ell-\alpha)\underline{x}\,\omega_{\ell-\alpha,\beta}(\underline{x})+2\beta\underline{x}\omega_{\ell-\alpha,\beta}(\underline{x})\, \omega_{1,0}(\underline{x})]\partial_{\underline{x}}^{\ell}\omega_{\alpha,-\beta}(\underline{x})\\
\medskip&=&(-1)^\ell \,\omega_{\ell-\alpha,\beta}(\underline{x}) [-(2(\alpha-\ell)\underline{x}-2\beta\underline{x}
\omega_{1,0}(\underline{x})]\,\partial_{\underline{x}}^{\ell}\omega_{\alpha,-\beta}(\underline{x})\\
\medskip&=&\Lambda(\underline{x}).
\end{array}
$$
As a result,
$$
K_{\ell+1,m}^{\alpha,-\beta}(\underline{x})=(-1)^{\ell+1} \omega_{\ell-\alpha+1,\beta}(\underline{x})\,
\partial_{\underline{x}}^{\ell+1}\omega_{\alpha,-\beta}(\underline{x}).
$$
\begin{prop}\label{orthogonalityRelation2} Let
$$
I_{\ell,t,p}^{\alpha,-\beta}=\displaystyle\int_{\mathbb{R}^m}\underline{x}^{\ell} K_{t,m}^{\alpha+p,\beta}(\underline{x})\, \omega_{\alpha,-\beta}(\underline{x})\, dV(\underline{x}).
$$
Whenever $2t<1-m-2\alpha$ we have the orthogonality relation
\begin{equation}\label{21}
I_{\ell,t,t}^{\alpha,-\beta}=0.
\end{equation}	
\end{prop}
\hskip-20pt\textbf{Proof.} Denote
$$
I_{\ell,t}=\displaystyle\int_{\mathbb{R}^m} \underline{x}^{\ell}\,\partial_{\underline{x}}^{t}(\omega_{\alpha+t,\beta}(\underline{x}))\,dV(\underline{x})).
$$
Using {Stokes's theorem, we obtain
$$
\begin{array}{lll}
\medskip&&\displaystyle\int_{\mathbb{R}^m}\underline{x}^{\ell} K_{t,m}^{\alpha+t,\beta}(\underline{x})\,\omega_{\alpha,-\beta}(\underline{x})dV(\underline{x})\\
\medskip&=&(-1)^t\,\displaystyle\int_{\mathbb{R}^m} \underline{x}^{\ell}\omega_{t-\alpha-t,\beta}(\underline{x}) \,\partial_{\underline{x}}^t(\omega_{\alpha+t,\beta}(\underline{x}))\,\omega_{\alpha,-\beta}(\underline{x})\,dV(\underline{x})\\
\medskip&=&(-1)^t\displaystyle\int_{\mathbb{R}^m} \underline{x}^{\ell}\,\partial_{\underline{x}}\partial_{\underline{x}}^{t-1}(\omega_{\alpha+t,\beta}(\underline{x}))\,dV(\underline{x}))\\
\medskip&=&(-1)^t\left[\displaystyle\int_{\partial\mathbb{R}^m}\underline{x}^{\ell}\partial_{\underline{x}}^{t-1}
(\omega_{\alpha+t,\beta}(\underline{x}))\,dV(\underline{x})\right.\\
\medskip&&\qquad\qquad\left.-\displaystyle\int_{\mathbb{R}^m}\partial_{\underline{x}}(\underline{x}^{\ell})
\partial_{\underline{x}}^{t-1}(\omega_{\alpha+t,\beta}(\underline{x}))\,dV(\underline{x})\right].
\end{array}
$$
Denote here-also
$$
I=\displaystyle\int_{\partial\mathbb{R}^m}\underline{x}^{\ell}
\partial_{\underline{x}}^{t-1}(\omega_{\alpha+t,\beta}(\underline{x}))\,dV(\underline{x})
$$
and
$$
II=\displaystyle\int_{\mathbb{R}^m}\partial_{\underline{x}}(\underline{x}^{\ell})\partial_{\underline{x}}^{t-1}
(\omega_{\alpha+t,\beta}(\underline{x}))\,dV(\underline{x}).
$$
The integral $I$ vanishes due to the assumption $2t<1-m-2\alpha$ and for the next, we have
$$
\begin{array}{lll}
\medskip\,II&=&\gamma_{l,m}\displaystyle\int_{\mathbb{R}^m}\underline{x}^{\ell-1}\partial_{\underline{x}}^{t-1}(\omega_{\alpha+t,\beta}(\underline{x}))\,dV(\underline{x})\\
\medskip&=& \gamma_{l,m} I_{\ell-1,t-1}.
\end{array}
$$
Hence, we obtain
$$
\begin{array}{lll}
\medskip\displaystyle\int_{\mathbb{R}^m}\underline{x}^{\ell} K_{t,m}^{\alpha+t,\beta}(\underline{x})\,\omega_{\alpha,-\beta}(\underline{x}) dV(\underline{x})&=&(-1)^{t+1}\gamma_{l,m}I_{\ell-1,t-1}\\
\medskip&=&(-1)^{t+1}\gamma_{l,m}[(-1)^{t}\gamma_{l-1,m}I_{\ell-2,t-2}]\\
&=&(-1)^{2t+1}\gamma_{l,m}\gamma_{l-1,m}I_{\ell-2,t-2}\\
&\vdots&\\
&=&C(m,\ell,t)\,I_0\\
&=&0
\end{array}
$$
where $C(m,\ell,t)$ is defined by (\ref{C(m,l,t)}).
\begin{defn}\label{GCGGJWaveletMother}
The generalized Clifford-Gauss-Gagenbauer-Jacobi Wavelet mother is defined by
$$
\psi_{\ell,m}^{\alpha,-\beta}(\underline{x})
=K_{\ell,m}^{\alpha+\ell,\beta}(\underline{x})
\omega_{\alpha,-\beta}(\underline{x})
=(-1)^\ell\partial_{\underline{x}}^{(\ell)}\,\omega_{\alpha+\ell,-\beta}(\underline{x}).
$$
\end{defn}
As for the previous class, we may prove that the wavelet $\psi_{\ell,m}^{\alpha,-\beta}(\underline{x})$ posses further vanishing moments as is shown in the next proposition.
\begin{prop} The following assertions are true.
\begin{enumerate}
\item Whenever $0<k<-m-\ell-2\alpha$ and $k<\ell$, we have
\begin{equation}\label{22}
\displaystyle\int_{\mathbb{R}^m} \underline{x}^k \psi_{\ell,m}^{\alpha,-\beta}(\underline{x}) dV(\underline{x})=0.
\end{equation}
\item The Clifford-Fourier transform of the generalized Clifford-Gauss-Gegenbauer-Jacobi wavelet is
\begin{equation}
\widehat{\psi_{\ell,m}^{\mu,\alpha}}(\underline{u})=(-i)^\ell\, \xi^\ell(2\pi)^{\frac{m}{2}}\rho^{1-\frac{m}{2}+\ell}\,\displaystyle\int_{0}^\infty\,\widetilde{\omega}_{\alpha,\beta}^{l,m}(r)\,J_{\frac{m}{2}-1}(r\rho)dr
\end{equation}	
where
$$
\widetilde{\omega}_{\alpha,\beta}^{l,m}(r)=(1+r^2)^{\alpha+\ell} r^{\frac{m}{2}} e^{-\beta r^2}.
$$
\end{enumerate}
\end{prop}
\begin{defn}
Let $a>0$ and $\underline{b}\in\mathbb{R}^m$. The copy of the generalized Clifford-Gauss-Gegenbauer-Jacobi wavelet mother at the scale $a$ and the position $\textbf{}$ is defined by
$$
_a^{\underline{b}}\psi_{\ell,m}^{\alpha,-\beta}(\underline{x})=a^{-\frac{m}{2}}\psi_{\ell,m}^{\alpha,-\beta}(\displaystyle\frac{\underline{x}-\underline{b}}{a}).
$$
The generalized Clifford-Gauss-Gegenbauer-Jacobi wavelet transform of a function $f\in L_2$ is defined by
$$
C_{a,\underline{b}}(f)=<_a^{\underline{b}}\psi_{\ell,m}^{\alpha,-\beta},f>.
$$
\end{defn}
\hskip-20pt The following result holds.
\begin{thm}
Let $\psi_{\ell,m}^{\mu,\alpha}$ be the wavelet defined in Definition \ref{GCGGJWaveletMother}. The following assertions hold.
\begin{enumerate}
\item
\begin{equation}\label{49}
\mathcal{A}_\psi=\displaystyle\frac{1}{\omega_m}\displaystyle\int_{\mathbb{R}^m}|\widehat{\psi_{\ell,m}^{\alpha,-\beta}}(\underline{x})|^2\displaystyle\frac{dV(\underline{x})}{|\underline{x}|^m}<+\infty.
\end{equation}
\item Any $L_2$ function $f$ may be reconstructed as
\begin{equation}\label{13}
f(x)=\displaystyle\frac{1}{\mathcal{A}_\psi}\displaystyle\int_{a>0}\displaystyle\int_{b\in\mathbb{R}^m}C_{a,\underline{b}}(f)\psi\left(\displaystyle\frac{\underline{x}-\underline{b}}{a}\right)\displaystyle\frac{da\,dV(\underline{b})}{a^{m+1}}.
\end{equation}
\end{enumerate}
\end{thm}
\hskip-20pt\textbf{Proof.}
\section{Conclusion}
In this paper new classes of monogenic orthogonal polynomials have been introduced relatively to different weights in the context of Clifford analysis. The new classes generalize the well known Jacobi and Gegenbauer polynomials. Such polynomial are proved to be good candidates to construct new wavelets in Clifford analysis. Fourier-Plancherel type results are generalized for the new classes of wavelets.

\end{document}